# Some facts about functionals of location and scatter

## R. M. Dudley[1],[*]

*Massachusetts Institute of Technology*

**Abstract:** Assumptions on a likelihood function, including a local Glivenko-Cantelli condition, imply the existence of M-estimators converging to an M-functional. Scatter matrix-valued estimators, defined on all empirical measures on $\mathbb{R}^d$ for $d \geq 2$, and equivariant under all, including singular, affine transformations, are shown to be constants times the sample covariance matrix. So, if weakly continuous, they must be identically 0. Results are stated on existence and differentiability of location and scatter functionals, defined on a weakly dense, weakly open set of laws, via elliptically symmetric t distributions on $\mathbb{R}^d$, following up on work of Kent, Tyler, and Dümbgen.

## 1. Introduction

In this paper a *law* will be a Borel probability measure on $\mathbb{R}^d$. Let $\mathcal{N}_d$ be the set of all $d \times d$ nonnegative definite symmetric matrices and $\mathcal{P}_d \subset \mathcal{N}_d$ the subset of strictly positive definite symmetric matrices. For $(\mu, \Sigma) \in \Theta = \mathbb{R}^d \times \mathcal{N}_d$, $\mu$ will be viewed as a location parameter and $\Sigma$ as a scatter parameter, extending the notions of mean vector and covariance matrix to arbitrarily heavy-tailed distributions. For $d \geq 2$, $\Theta$ may be taken to be $\mathcal{P}_d$ or $\mathbb{R}^d \times \mathcal{P}_d$.

For a law $P$ on $\mathbb{R}^d$, let $X_1, X_2, \ldots$ be i.i.d. $(P)$ and let $P_n$ be the empirical measure $n^{-1} \sum_{j=1}^{n} \delta_{X_j}$ where $\delta_x(A) := 1_A(x)$ for any point $x$ and set $A$. A class $\mathcal{F} \subset \mathcal{L}^1(\mathbb{R}^d, P)$ is called a *Glivenko-Cantelli* class for $P$ if

$$(1) \qquad \sup\{|\int f d(P_n - P)| : f \in \mathcal{F}\} \to 0$$

almost surely as $n \to \infty$ (if the supremum is measurable, as it will be in all cases considered in this paper). Talagrand [20, 21] characterized such classes. A class $\mathcal{F}$ of Borel measurable functions on $\mathbb{R}^d$ is called a *universal* Glivenko-Cantelli class if it is a Glivenko-Cantelli class for *all* laws $P$ on $\mathbb{R}^d$, and a *uniform* Glivenko-Cantelli class if the convergence in (1) is uniform over all laws $P$. Rather general sufficient conditions for the universal Glivenko-Cantelli property and a characterization up to measurability of the uniform property have been given [7].

Let $\rho : (x, \theta) \mapsto \rho(x, \theta) \in \mathbb{R}$ defined for $x \in \mathbb{R}^d$ and $\theta \in \Theta$, Borel measurable in $x$ and lower semicontinuous in $\theta$, i.e. $\rho(x, \theta) \leq \liminf_{\phi \to \theta} \rho(x, \phi)$ for all $x$ and $\theta$. For a law $Q$, let $Q\rho(\phi) := \int \rho(x, \phi) dQ(x)$ if the integral is defined (not $\infty - \infty$), as it always will be if $Q = P_n$. An *M-estimate* of $\theta$ for a given $n$ and $P_n$ will be a $\hat{\theta}_n$ such that $P_n \rho(\theta)$ is minimized at $\theta = \hat{\theta}_n$, if it exists and is unique. A measurable function, not necessarily defined a.s., whose values are M-estimates is called an M-estimator. An *M-limit* $\theta_0 = \theta_0(P) = \theta_0(P, \rho)$ (with respect to $\rho$) will mean a point of $\Theta$ such for every open neighborhood $U$ of $\theta_0$, as $n \to \infty$,

$$(2) \qquad \Pr\Big\{\inf\{P_n\rho(\theta) : \ \theta \notin U\} \leq \inf\{P_n\rho(\phi) : \ \phi \in U\}\Big\} \to 0,$$

where the given probabilities are assumed to be defined. Then if M-estimators exist (with probability $\to 1$ as $n \to \infty$), they must converge in probability to $\theta_0(P)$. An M-limit $\theta_0 = \theta_0(P)$ with respect to $\rho$ will be called *definite* iff for every neighborhood $U$ of $\theta_0$ there is an $\varepsilon > 0$ such that the outer probability

$$(3) \qquad (P^n)^* \{\inf\{P_n\rho(\theta) : \ \theta \notin U\} \leq \varepsilon + \inf\{P_n\rho(\phi) : \ \phi \in U\}\} \to 0$$

as $n \to \infty$.

[1]Room 2-245, MIT, Cambridge, MA 02139-4307, USA, e-mail: rmd@math.mit.edu
[*]Partially supported by NSF Grants DMS-0103821 and DMS-0504909.

For a law $P$ on $\mathbb{R}^d$ and a given $\rho(\cdot, \cdot)$, a $\theta_1 = \theta_1(P)$ is called the *M-functional* of $P$ for $\rho$ if and only if there exists a measurable function $a(x)$, called an *adjustment function*, such that for $h(x, \theta) = \rho(x, \theta) - a(x)$, $Ph(\theta)$ is defined and satisfies $-\infty < Ph(\theta) \leq +\infty$ for all $\theta \in \Theta$, and is minimized uniquely at $\theta = \theta_1(P)$, e.g. Huber [13]. As Huber showed, $\theta_1(P)$ doesn't depend on the choice of $a(\cdot)$. Clearly, an M-estimate $\hat{\theta}_n$ is the M-functional $\theta_1(P_n)$ if either exists.

A lower semicontinuous function $f$ from $\Theta$ into $(-\infty, +\infty]$ will be called *uniminimal* iff it has a unique relative minimum at a point $\theta_0$ and for all $t \in \mathbb{R}$, $\{\theta \in \Theta : \ f(\theta) \leq t\}$ is connected. For a differentiable function $f$, recall that a *critical point* of $f$ is a point where the gradient of $f$ is 0.

**Examples.** On $\Theta = \mathbb{R}$ let $f(x) = -(1 - x^2)^2$. Then $f$ has a unique relative minimum at $x = 0$, but no absolute minimum. It has two other critical points which are relative maxima. For $t < 0$ the set where $f \leq t$ is not connected.

If $f$ is a strictly convex function on $\mathbb{R}^d$ attaining its minimum, then $f$ is uniminimal, as is $\theta \mapsto f(x - \theta)$ for any $x$. So is $\theta \mapsto \int f(x - \theta) - f(x) dP(x)$ if it's defined and finite and attains its minimum for a law $P$, as will be true e.g. if $f(x) = |x|^2$ and $\int |x| dP(x) < \infty$, or for all $P$ if $f$ is also Lipschitz, e.g. $f(x) = \sqrt{1 + |x|^2}$.

I have not found the notion here called "uniminimal" in the literature. Similar but more complex assumptions occur in some work on sufficient conditions for minimaxity in game theory, e.g. [11]. Thus, I claim no originality for the following easily proved fact.

**Proposition 1.** *Let* $(\Theta, d)$ *be a locally compact metric space. If* $f$ *is uniminimal on* $(\Theta, d)$, *then (a)* $f$ *attains its absolute minimum at its unique relative minimum* $\theta_0$, *and (b) For every neighborhood* $U$ *of* $\theta_0$ *there is an* $\varepsilon > 0$ *such that* $f(\theta) \geq f(\theta_0) + \varepsilon$ *for all* $\theta \notin U$.

*Proof.* Clearly (b) implies (a). To prove (b), suppose that for some or equivalently all small enough $\delta > 0$ and all $n = 1, 2, \ldots$, there are $\theta_n \in \Theta$ with $d(\theta_n, \theta_0) \geq \delta$ and $f(\theta_n) \leq f(\theta_0) + 1/n$. By connectedness, we can take $d(\theta_n, \theta_0) = \delta$ for all $n$. Then for $\delta > 0$ small enough, $\{\theta : \ d(\theta, \theta_0) \leq \delta\}$ is compact and there is a converging subsequence $\theta_{n(k)} \to \theta_\delta$ with $d(\theta_\delta, \theta_0) = \delta$ and $f(\theta_\delta) \leq f(\theta_0)$ by lower semicontinuity. Letting $\delta \downarrow 0$ we get a contradiction to the fact that $\theta_0$ is a unique relative minimum. $\square$

**Theorem 2.** *Let* $(\Theta, d)$ *be a connected locally compact metric space and* $(X, \mathcal{B}, P)$ *a probability space. Let* $h : \ X \times \Theta \mapsto \mathbb{R}$ *where for each* $\theta \in \Theta$, $h(\cdot, \theta)$ *is measurable. Assume that:*

(i) $\theta \mapsto Ph(\theta) \in (-\infty, +\infty]$ *is well-defined and uniminimal on* $\Theta$, *with minimum at* $\theta_0$;

(ii) *Outside an event* $A_n$ *whose probability converges to 0 as* $n \to \infty$, $P_n h(\cdot)$ *is uniminimal on* $\Theta$;

(iii) *For some neighborhood* $U$ *of* $\theta_0$, $\{h(\cdot, \theta) : \theta \in U\}$ *is a Glivenko-Cantelli class for* $P$.

Then $\theta_0$ *is the definite M-limit for* $P$ *and the M-functional* $\theta_1(P)$.

**Remark.** Glivenko-Cantelli conditions on log likelihoods (and their partial derivatives through order 2) for parameters in bounded neighborhoods have been assumed in other work, e.g. [17] and [8].

*Proof.* That $\theta_0$ is an M-functional for $P$ follows from (i) and Proposition 1. By (iii), take $\delta > 0$ small enough so that $\{h(\cdot, \theta) : d(\theta, \theta_0) < \delta\}$ is a Glivenko-Cantelli class for $P$. By (i) and Proposition 1, take $\varepsilon > 0$ such that $Ph(\theta) > Ph(\theta_0) + 3\varepsilon$ whenever $d(\theta, \theta_0) > \delta/2$. Outside some events $A_n$ whose probability converges to 0 as $n \to \infty$, we have $P_n h(\theta_0) < Ph(\theta_0) + \varepsilon$ and $P_n h(\theta) > Ph(\theta_0) + 2\varepsilon$ for all $\theta$ with $\delta/2 < d(\theta, \theta_0) < \delta$. Then by (ii), also with probability converging to 1, $P_n h(\theta) > P_n h(\theta_0) + \varepsilon$ for all $\theta$ with $d(\theta, \theta_0) > \delta/2$, proving (3) and the theorem. $\quad\square$

A class $\mathcal{C}$ of subsets of a set $X$ is called a *VC (Vapnik-Chervonenkis) class* if for some $k < \infty$, for every subset $A$ of $X$ with $k$ elements, there is some $B \subset A$ with $B \neq C \cap A$ for all $C \in \mathcal{C}$, e.g. [4, Chapter 4]. A class $\mathcal{F}$ of real-valued functions on $X$ is called a *VC major class* iff $\{\{x \in X : f(x) > t\} : f \in \mathcal{F}, \ t \in \mathbb{R}\}$ is a VC class of sets (e.g. [4, Section 4.7]). In the following, local compactness is stronger than needed but holds for the parameter spaces being considered.

**Theorem 3.** *Let* $h(x, \theta)$ *be continuous in* $\theta \in \Theta$ *for each* $x$ *and measurable in* $x$ *for each* $\theta$ *where* $\Theta$ *is a locally compact separable metric space. Let* $h(\cdot, \cdot)$ *be uniformly bounded and let* $\mathcal{F} := \{h(\cdot, \theta) : \theta \in \Theta\}$ *be a VC major class of functions. Then* $\mathcal{F}$ *is a uniform, thus universal, Glivenko-Cantelli class.*

*Proof.* Theorem 6 of [7] applies: sufficient bounds for the Koltchinskii-Pollard entropy of uniformly bounded VC major classes of functions are given in [3, Theorem 2.1(a), Corollary 5.8], and sufficient measurability of the class $\mathcal{F}$ follows from the continuity in $\theta$ and the assumptions on $\Theta$. $\quad\square$

For the $t$ location-scatter functionals in Sections 4 and 5, the notions of VC major class, and local Glivenko-Cantelli class as in Theorem 2(iii), will be applicable. But as shown by Kent, Tyler and Vardi [16], to be recalled after Theorem 12(iii), some parts of the development work only for $t$ functionals, rather than for functions $\rho$ satisfying general properties such as convexity.

## 2. Equivariance for location and scatter

Notions of "location" and "scale" or multidimensional "scatter" functional will be defined along with equivariance, as follows.

**Definitions.** Let $Q \mapsto \mu(Q) \in \mathbb{R}^d$, resp. $\Sigma(Q) \in \mathcal{N}_d$, be a functional defined on a set $\mathcal{D}$ of laws $Q$ on $\mathbb{R}^d$. Then $\mu$ (resp. $\Sigma$) is called an *affinely equivariant location* (resp. *scatter*) *functional* iff for any nonsingular $d \times d$ matrix $A$ and $v \in \mathbb{R}^d$, with $f(x) := Ax + v$, and any law $Q \in \mathcal{D}$, the image measure $P := Q \circ f^{-1} \in \mathcal{D}$ also, with $\mu(P) = A\mu(Q) + v$ or, respectively, $\Sigma(P) = A\Sigma(Q)A'$. For $d = 1$, $\sigma(\cdot)$ with $0 \leq \sigma < \infty$ will be called an *affinely equivariant scale functional* iff $\sigma^2$ satisfies the definition of affinely equivariant scatter functional.

Well-known examples of affinely equivariant location and scale functionals (for $d = 1$), defined for all laws, are the median and MAD (median absolute deviation), where for a real random variable $X$ with median $m$, the MAD of $X$ or its distribution is defined as the median of $|X - m|$.

Call a location functional $\mu(\cdot)$ or a scatter functional $\Sigma(\cdot)$ *singularly affine equivariant* if in the definition of affine equivariance $A$ can be any matrix, possibly singular. If a functional is defined on all laws, affinely equivariant, and weakly continuous, then it must be singularly affine equivariant. The classical sample mean and covariance are defined for all $P_n$ and singularly affine equivariant. It turns out that in dimension $d \geq 2$, there are essentially no other singularly affine equivariant location or scatter functionals defined for all $P_n$, and so weak continuity at all laws is not possible. First the known fact for location will be recalled, then an at least partially known fact for scatter will be stated and proved.

Let $X$ be a $d \times n$ data matrix whose $j$th column is $X_j \in \mathbb{R}^d$. Let $X^i$ be the $i$th row of $X$. Let $\mathbf{1}_n$ be the $n \times 1$ vector with all components 1. Let $\overline{X} = \int x dP_n$ be the sample mean vector in $\mathbb{R}^d$, so that $X - \overline{X}\mathbf{1}'_n$ is the centered data matrix. Note that $P_n$, and thus $\overline{X}$ and $\Sigma(X)$, are preserved by any permutation of the columns of $X$.

**Theorem 4.** (a) *If $\mu(\cdot)$ is a singularly affine equivariant location functional defined for all $P_n$ on $\mathbb{R}^d$ for $d \geq 2$ and a fixed $n$, then $\mu(P_n) \equiv \overline{X}$.*
(b) *If in addition $\mu(\cdot)$ is defined for all $n$ and all $P_n$ on $\mathbb{R}^d$, then as $n$ varies, $\mu(\cdot)$ is not weakly continuous. Thus, there is no affinely equivariant, weakly continuous location functional defined on all laws on $\mathbb{R}^d$ for $d \geq 2$.*

*Proof.* Part (a) follows from work of Obenchain [18, Lemma 1] and permutation invariance, as noted e.g. by Rousseeuw [19]. Then (b) follows directly, for $x_1 = n$, $x_2 = \cdots = x_n = 0$, $n \to \infty$.                                                    □

Next is a related fact about scatter functionals. Davies [1, p. 1879] made a statement closely related to part (b), strong but not quite in the same generality, and very briefly suggested a proof. I don't know a reference for part (a), or an explicit one for (b), so a proof will be given.

**Theorem 5.** (a) *Let $\Sigma(\cdot)$ be a singularly affine equivariant scatter functional defined on all empirical measures $P_n$ on $\mathbb{R}^d$ for $d \geq 2$ and some fixed $n \geq 2$. Write $\Sigma(X) := \Sigma(P_n)$. Then there is a constant $c_n \geq 0$, depending on $\Sigma(\cdot)$, such that for any $X$, $\Sigma(X - \overline{X}\mathbf{1}'_n) = c_n(X - \overline{X}\mathbf{1}'_n)(X - \overline{X}\mathbf{1}'_n)'$. In other words, applied to centered data matrices, $\Sigma$ is proportional to the sample covariance matrix.*
(b) *If $\Sigma(\cdot)$ is an affinely equivariant scatter functional defined for all $n$ and $P_n$ on $\mathbb{R}^d$ for $d \geq 2$, weakly continuous as a function of $P_n$, then $\Sigma \equiv 0$.*

*Proof.* (a) We have $\Sigma(BX) = B\Sigma(X)B'$ for any $d \times d$ matrix $B$. For any $U, V \in \mathbb{R}^n$ let $X^1 = U'$, $X^2 = V'$, and $(U, V) := \Sigma_{12}(X)$. Then $(\cdot, \cdot)$ is well-defined, letting $B_{11} = B_{22} = 1$ and $B_{ij} = 0$ otherwise. It will be shown that $(\cdot, \cdot)$ is a semi-inner product. We have $(U, V) \equiv (V, U)$ via $B$ with $B_{12} = B_{21} = 1$ and $B_{ij} = 0$ otherwise, since $\Sigma$ is symmetric. For $B_{11} = B_{21} = 1$ and $B_{ij} = 0$ otherwise we get for any $U \in \mathbb{R}^n$ that

$$(4) \qquad (U, U) = \Sigma_{12}(BX) = (B\Sigma(X)B')_{12} = \Sigma_{11}(X) \geq 0.$$

For constants $a$ and $b$, $(aU, bV) \equiv ab(U, V)$ follows for $B_{11} = a$, $B_{22} = b$, and $B_{ij} = 0$ otherwise. It remains to prove biadditivity $(U, V + W) \equiv (U, V) + (U, W)$. For $d \geq 3$ this is easy, letting $X^3 = W$, $B_{11} = B_{22} = B_{23} = 1$, and $B_{ij} = 0$ otherwise. For $d = 2$, we first get $(U + V, V) = (U, V) + (V, V)$ from $B = \left(\begin{smallmatrix} 1 & 1 \\ 0 & 1 \end{smallmatrix}\right)$. Symmetrically, $(U, U + V) = (U, U) + (U, V)$. Then from $B = \left(\begin{smallmatrix} 1 & 1 \\ 1 & 1 \end{smallmatrix}\right)$ we get

$$(5) \qquad (U + V, U + V) = (U, U) + 2(U, V) + (V, V).$$

Letting $\|W\|^2 := (W, W)$ for any $W \in \mathbb{R}^n$ we get the parallelogram law $\|U + V\|^2 + \|U - V\|^2 \equiv 2\|U\|^2 + 2\|V\|^2$. (But $\|\cdot\|$ has not yet been shown to be a norm.) Applying this repeatedly we get for any $W, Y$, and $Z \in \mathbb{R}^n$ that

$$\|W + Y + Z\|^2 - \|W - Y - Z\|^2 = \|W + Y\|^2 - \|W - Y\|^2 + \|W + Z\|^2 - \|W - Z\|^2,$$

letting first $U = W + Y$, $V = Z$, then $U = W - Z$, $V = Y$, then $U = W$, $V = Z$, and lastly $U = W$, $V = Y$. Applying (5) gives $(W, Y + Z) \equiv (W, Y) + (W, Z)$, the desired biadditivity. So $(\cdot, \cdot)$ is indeed a semi-inner product, i.e. there is a $C(n) \in \mathcal{N}_n$ such that $(U, V) \equiv U' C(n) V$. By permutation invariance, there are numbers $a_n \geq 0$ and $b_n$ such that $C(n)_{ii} = a_n$ for all $i = 1, \ldots, n$ and $C(n)_{ij} = b_n$ for all $i \neq j$.

Let $c_n := a_n - b_n$ and let $e_i \in \mathbb{R}^n$ be the $i$th standard unit vector. For each $y \in \mathbb{R}^n$ let $y = \sum_{i=1}^n y_i e_i$ and $\overline{y} := \frac{1}{n} \sum_{i=1}^n y_i$. Then for any $z \in \mathbb{R}^n$,

$$(y - \overline{y}\mathbf{1}_n, z - \overline{z}\mathbf{1}_n) = \sum_{i,j=1}^n C(n)_{ij}(y_i - \overline{y})(z_j - \overline{z}) = c_n(y - \overline{y}\mathbf{1}_n)'(z - \overline{z}\mathbf{1}_n).$$

For $1 \leq j \leq k \leq d$, let $B_{ir} := \delta_{r\pi(i)}$ for a function $\pi$ from $\{1, 2, \ldots, d\}$ into itself with $\pi(1) = j$ and $\pi(2) = k$. Then $(BX)^1 = X^j$ and $(BX)^2 = X^k$. Thus $(X^j, X^k) = \Sigma_{12}(BX) = \Sigma_{jk}(X)$, recalling (4) if $j = k$.

Let $\overline{X} \in \mathbb{R}^d$ have $i$th component $\overline{X}^i$. Then

$$\Sigma_{jk}(X - \overline{X}\mathbf{1}'_n) = (X^j - \overline{X}^j\mathbf{1}_n, X^k - \overline{X}^k\mathbf{1}_n) = c_n(X^j - \overline{X}^j\mathbf{1}_n)'(X^k - \overline{X}^k\mathbf{1}_n),$$

where $c_n \geq 0$ is seen when $j = k$ and the coefficient of $c_n$ is strictly positive, as it can be since $n \geq 2$. Thus part (a) is proved.

For part (b), consider empirical measures $P_n = P_{mn}$, so that each $X_j$ in $P_n$ is repeated $m$ times in $P_{mn}$. Since the $\overline{X}$'s and $\Sigma$s for $P_n$ and $P_{mn}$ must be the same, we get that $c_{mn} = c_n/m$ which likewise equals $c_m/n$. Thus there is a constant $c_1$ such that $c_n = c_1/n$ for all $n$.

Let $X_{11} := -X_{12} := \sqrt{n}$, let $X_{ij} = 0$ for all other $i, j$ and let $n \to \infty$. Then $\overline{X} \equiv 0$, $P_n \to \delta_0$ weakly, and $\Sigma(\delta_0)$ is the 0 matrix by singular affine equivariance with $B = 0$, but $\Sigma(P_n)$ don't converge to 0 unless $c_1 = 0$ and so $c_n = 0$ for all $n$, proving (b). □

So, for $d \geq 2$, affinely equivariant location and non-zero scatter functionals, weakly continuous on their domains, can't be defined on all laws. They can be defined on weakly dense and open domains, as will be seen in Theorem 12, on which they can have good differentiability properties, as seen in Section 5.

## 3. Multivariate scatter

This section treats pure scatter in $\mathbb{R}^d$, with $\Theta = \mathcal{P}_d$. Results of Kent and Tyler [15] for finite samples, to be recalled, are extended to general laws on $\mathbb{R}^d$ in [6, Section 3].

For $A \in \mathcal{P}_d$ and a function $\rho$ from $[0, \infty)$ into itself, consider the function

$$(6) \qquad L(y, A) \; := \; \frac{1}{2} \log \det A + \rho(y' A^{-1} y), \quad y \in \mathbb{R}^d.$$

For adjustment, let

$$(7) \qquad h(y, A) \; := \; L(y, A) - L(y, I)$$

where $I$ is the identity matrix. Then

$$(8) \qquad Qh(A) \; = \; \frac{1}{2} \log \det A + \int \rho(y' A^{-1} y) - \rho(y' y) \, dQ(y)$$

if the integral is defined. We have the following, shown for $Q = Q_n$ an empirical measure in [15, (1.3)] and for general $Q$ in [6, Section 3]. Here (9) is a redescending condition. A symmetric $d \times d$ matrix $A$ will be parameterized by the entries $A_{ij}$ for $1 \leq i \leq j \leq d$. Thus in taking a partial derivative of a function $f(A)$ with respect to an entry $A_{ij}$, $A_{ji} \equiv A_{ij}$ will vary while $A_{kl}$ will remain fixed except for $(k, l) = (i, j)$ or $(j, i)$.

**Proposition 6.** *Let $\rho$ be continuous from $[0, \infty)$ into itself and have a bounded continuous derivative, where $\rho'(0) \; := \; \rho'(0+) \; := \; \lim_{x \downarrow 0} [\rho(x) - \rho(0)]/x$. Let $0 \leq u(x) \; := \; 2\rho'(x)$ for $x \geq 0$. Assume that*

$$(9) \qquad \sup_{0 \leq x < \infty} x u(x) < \infty.$$

*Then for each law $Q$ on $\mathbb{R}^d$, $Qh$ in (8) is well defined and is a $C^1$ function of the entries of $A$. Here $Qh$ has a critical point at $A = B$ if and only if*

$$(10) \qquad B \; = \; \int u(y' B^{-1} y) y y' dQ(y).$$

The following, proved in [6, Section 3], extends to any law $Q$ the uniqueness part of [15, Theorem 2.2].

**Proposition 7.** *Under the hypotheses of Proposition 6, if in addition $u(\cdot)$ is non-increasing and $s \mapsto su(s)$ is strictly increasing on $[0, \infty)$, then for any law $Q$ on $\mathbb{R}^d$, $Qh$ has at most one critical point $A \in \mathcal{P}_d$.*

A sufficient condition for existence of a pure scatter M-functional $A(Q)$ will include the following assumption from [15, (2.4)]. Given a function $u(\cdot)$ as in Proposition 7, let $a_0 \; := \; a_0(u(\cdot)) \; := \; \sup_{s > 0} su(s)$. Since $s \mapsto su(s)$ is increasing, it follows that

$$(11) \qquad su(s) \uparrow a_0 \quad \text{as} \quad s \uparrow + \infty.$$

Kent and Tyler [15] gave the following condition for empirical measures.

**Definition.** Given $a_0 \; := \; a(0) > 0$, let $\mathcal{U}_{d, a(0)}$ denote the set of all laws $Q$ on $\mathbb{R}^d$ such that for every proper linear subspace $H$ of $\mathbb{R}^d$, of dimension $q \leq d - 1$, we have $Q(H) < 1 - (d - q)/a_0$.

Note that $\mathcal{U}_{d, a(0)}$ is weakly open and dense and contains all laws with densities. If $Q \in \mathcal{U}_{d, a(0)}$, then $Q(\{0\}) < 1 - (d/a_0)$, which is impossible if $a_0 \leq d$. So in the next theorem we assume $a_0 > d$. In part (b), the existence of a unique $B(Q_n)$ minimizing $Q_n h$ for an empirical $Q_n \in \mathcal{U}_{d, a(0)}$ was proved in [15, Theorems 2.1 and 2.2]. For a general $Q \in \mathcal{U}_{d, a(0)}$ it's proved in [6, Section 3]; one lemma useful in the proof is proved here.

**Theorem 8.** *Under the assumptions of Propositions* 6 *and* 7, *for* $a(0) = a_0$ *as in* (11),

(a) *If* $Q \notin \mathcal{U}_{d,a(0)}$, *then* $Qh$ *has no critical points.*

(b) *If* $a_0 > d$ *and* $Q \in \mathcal{U}_{d,a(0)}$, *then* $Qh$ *attains its minimum at a unique* $B = B(Q) \in \mathcal{P}_d$ *and has no other critical points.*

A proof of the theorem uses a fact about probabilities of proper subspaces or hyperplanes. A related statement is Lemma 5.1 of Dümbgen and Tyler [10].

**Lemma 9.** *Let* $V$ *be a real vector space with a* $\sigma$-*algebra* $\mathcal{B}$ *for which all finite-dimensional hyperplanes* $H = x + T := \{x + u : u \in T\}$ *for finite-dimensional vector subspaces* $T$ *are measurable. Let* $Q$ *be a probability measure on* $\mathcal{B}$ *and let* $\mathcal{H}_j$ *be the collection of all* $j$-*dimensional hyperplanes in* $V$. *Then for each* $j = 0, 1, 2, \ldots$, *for any infinite sequence* $\{C_i\}$ *of distinct hyperplanes in* $\mathcal{H}_j$ *such that* $Q(C_i)$ *converges, its limit must be* $Q(F)$ *for some hyperplane* $F$ *of dimension less than* $j$ *such that* $F \subset C_i$ *for infinitely many* $i$. *In particular,* $Q(C_i)$ *cannot be strictly increasing. The same is true for vector subspaces in place of hyperplanes.*

*Proof.* Hyperplanes of dimension 0 are singletons $\{x\}$. The empty set $\emptyset$ will be considered as a hyperplane of dimension $-1$. Let $W_{-1} := \emptyset$. *Claim 1*: For each $j = 0, 1, \ldots$, there exists a finite or countable sequence $\{V_{ji}\} \subset \mathcal{H}_j$ such that for $W_j := W_{j-1} \cup \bigcup_i V_{ji}$, $Q(V \setminus W_j) = 0$ for all $V \in \mathcal{H}_j$. Let $V_{0i} = \{x_i\}$ for some unique $i$ if and only if $Q(\{x_i\}) > 0$. The set of such $x_i$ is clearly countable. Let $W_0 := \cup_i V_{0i} = \{x \in V : Q(\{x\}) > 0\}$. Clearly, for any $x \in V$, $Q(\{x\} \setminus W_0) = 0$. Recursively, for $j \geq 1$, assuming $W_{j-1}$ has the given properties, suppose for $r = 1, 2$, $H_r \in \mathcal{H}_j$ and $Q(H_r \setminus W_{j-1}) > 0$. If $H_1 \neq H_2$, then $H_1 \cap H_2$ is a hyperplane of dimension at most $j-1$, so $Q(H_1 \cap H_2 \setminus W_{j-1}) = 0$ and the sets $H_r \setminus W_{j-1}$ are disjoint up to sets with $Q = 0$. Thus there are at most countably many different $H_r \in \mathcal{H}_j$ with $Q(H_r \setminus W_{j-1}) > 0$. Let $V_{jr} := H_r$ for such $H_r$ and set $W_j := W_{j-1} \cup \bigcup_r V_{jr}$. It's then clear that for any $H \in \mathcal{H}_j$, $Q(H \setminus W_j) = 0$, so the recursion can continue and Claim 1 is proved.

*Claim 2* is that if $C$ is any hyperplane of dimension $j$ or larger, and $s = 0, 1, \ldots, j$, then for each $r$, either $C \supset V_{sr}$ or $Q(C \cap (V_{sr} \setminus W_{s-1})) = 0$. If $C$ doesn't include $V_{sr}$, then $C \cap V_{sr}$ is a hyperplane of dimension $\leq s - 1$, and so included in $W_{s-1}$ up to a set with $Q = 0$, so Claim 2 follows.

Now, given distinct $C_i \in \mathcal{H}_j$ with $Q(C_i)$ converging, let $B$ be a hyperplane of largest possible dimension $b$ included in $C_i$ for infinitely many $i$. Then $b < j$. Taking a subsequence, we can assume that $B \subset C_i$ for all $i$. *Claim 3* is that then $Q(C_i \setminus B) \to 0$ as $i \to \infty$. For any $s = 0, 1, \ldots, j - 1$, and each $r$, by Claim 2, if $C_i \supset V_{sr}$ for infinitely many $i$, then $V_{sr} \subset B$, since otherwise $C_i$ includes the smallest hyperplane including $V_{sr}$ and $B$, which has dimension larger than $b$, a contradiction. So $\lim_{i \to \infty} Q((C_i \setminus B) \cap (V_{sr} \setminus W_{s-1})) = 0$ for each $s < j$ and $r$. It follows by induction on $s$ that $Q(C_i \cap W_s \setminus B) \to 0$ as $i \to \infty$ for $s = 0, 1, \ldots, j - 1$.

By the proof of Claim 1, the sets $C_i \setminus W_{j-1}$ are disjoint up to sets with $Q = 0$, so Claim 3 follows, and so the statement of the lemma for hyperplanes. The proof for vector subspaces is parallel and easier. The fact that $Q(C_i)$ cannot be strictly increasing then clearly follows, as a subsequence would also be strictly increasing. So the lemma is proved. $\square$

Dümbgen and Tyler [10], Lemma 5.1 show that $\sup\{Q(V) : V \in \mathcal{H}_j\}$ is attained for each $Q$ and $j$ and is weakly upper semicontinuous in $Q$.

### 4. Location and scatter $t$ functionals

As Kent and Tyler [15, Section 3] and Kent, Tyler and Vardi [16] showed, $(t)$ location-scatter estimation in $\mathbb{R}^d$ can be reduced to pure scatter estimation in $\mathbb{R}^{d+1}$, beginning with the following.

**Proposition 10.** (i) *For any $d = 1, 2, \ldots$, there is a 1-1 correspondence, $C^\infty$ in either direction, between matrices $A \in \mathcal{P}_{d+1}$ and triples $(\Sigma, \mu, \gamma)$ where $\Sigma \in \mathcal{P}_d$, $\mu \in \mathbb{R}^d$, and $\gamma > 0$, given by*

$$(12) \qquad A = A(\Sigma, \mu, \gamma) = \gamma \begin{bmatrix} \Sigma + \mu\mu' & \mu \\ \mu' & 1 \end{bmatrix}.$$

*The same holds for $A \in \mathcal{P}_{d+1}$ with $\gamma = A_{d+1, d+1} = 1$ and pairs $(\mu, \Sigma) \in \mathbb{R}^d \times \mathcal{P}^d$.*
(ii) *If (12) holds, then for any $y \in \mathbb{R}^d$ (a column vector),*

$$(13) \qquad (y', 1) A^{-1} (y', 1)' \;=\; \gamma^{-1} \left( 1 + (y - \mu)' \Sigma^{-1} (y - \mu) \right).$$

For M-estimation of location and scatter in $\mathbb{R}^d$, we will have a function $\rho : [0, \infty) \mapsto [0, \infty)$ as in the previous section. The parameter space is now the set of pairs $(\mu, \Sigma)$ for $\mu \in \mathbb{R}^d$ and $\Sigma \in \mathcal{P}_d$, and we have a multivariate $\rho$ function

$$\rho(y, (\mu, \Sigma)) \;:=\; \frac{1}{2} \log \det \Sigma + \rho((y - \mu)' \Sigma^{-1} (y - \mu)).$$

For any $\mu \in \mathbb{R}^d$ and $\Sigma \in \mathcal{P}_d$ let $A_0 := A_0(\mu, \Sigma) := A(\Sigma, \mu, 1) \in \mathcal{P}_{d+1}$ by (12) with $\gamma = 1$, noting that $\det A_0 = \det \Sigma$. Now $\rho$ can be adjusted, in light of (9) and (13), by defining

$$(14) \qquad h(y, (\mu, \Sigma)) \;:=\; \rho(y, (\mu, \Sigma)) - \rho(y'y).$$

Laws $P$ on $\mathbb{R}^d$ correspond to laws $Q := P \circ T_1^{-1}$ on $\mathbb{R}^{d+1}$ concentrated in $\{y : y_{d+1} = 1\}$, where $T_1(y) := (y', 1)' \in \mathbb{R}^{d+1}$, $y \in \mathbb{R}^d$. We will need a hypothesis on $P$ corresponding to $Q \in \mathcal{U}_{d+1, a(0)}$. Kent and Tyler [15] gave these conditions for empirical measures.

**Definition.** For any $a_0 > 0$ let $\mathcal{V}_{d, a(0)}$ be the set of all laws $P$ on $\mathbb{R}^d$ such that $P(J) < 1 - (d - q)/a_0$ for every affine hyperplane $J$ of dimension $q < d$.

The next fact is rather easy to prove. Here $a > d + 1$ avoids the contradictory $Q(\{0\}) < 0$.

**Proposition 11.** *If $P$ is a law on $\mathbb{R}^d$, $a > d + 1$, and $Q := P \circ T_1^{-1}$ on $\mathbb{R}^{d+1}$, then $P \in \mathcal{V}_{d, a}$ if and only if $Q \in \mathcal{U}_{d+1, a}$.*

A family of $\rho$ functions for which $\gamma = 1$ automatically, as noted by Kent and Tyler [15, (1.5), (1.6), Section 4], is given by elliptically symmetric multivariate $t$ densities with $\nu$ degrees of freedom as follows: for $0 < \nu < \infty$ and $0 \le s < \infty$ let

$$(15) \qquad \rho_\nu(s) \;:=\; \rho_{\nu, d}(s) \;:=\; \frac{\nu + d}{2} \log \left( \frac{\nu + s}{\nu} \right).$$

For this $\rho$, $u$ is $u_\nu(s) := u_{\nu, d}(s) := (\nu + d)/(\nu + s)$, which is decreasing, and $s \mapsto s u_{\nu, d}(s)$ is strictly increasing and bounded, i.e. (9) holds, with supremum and limit at $+\infty$ equal to $a_{0, \nu} := a_0(u_\nu(\cdot)) = \nu + d$.

The following fact is in part given by Kent and Tyler [15] and further by Kent, Tyler and Vardi [16], for empirical measures; equation (16) was not found explicitly in either. Here a proof will be given for any $P \in \mathcal{V}_{d, \nu+d}$, assuming Theorem 8 and Propositions 6 and 10.

**Theorem 12.** *For any $d = 1, 2, \ldots,$ $\nu > 1$, law $P$ on $\mathbb{R}^d$, and $Q = P \circ T_1^{-1}$ on $\mathbb{R}^{d+1}$, letting $\nu' := \nu - 1$, assuming $P \in \mathcal{V}_{d,\nu+d}$ in parts (a) through (e),*

(a) *For $A \in \mathcal{P}_{d+1}$, $A \mapsto Qh(A)$ defined by (8) for $\rho = \rho_{\nu',d+1}$ has a unique critical point $A(\nu') := A_{\nu'}(Q)$ which is an absolute minimum;*

(b) $A(\nu')_{d+1,d+1} = \int u_{\nu',d+1}(y'A(\nu')^{-1}y)dQ(y) = 1$;

(c) *For any $\mu \in \mathbb{R}^d$ and $\Sigma \in \mathcal{P}_d$ let $A = A(\Sigma, \mu, 1) \in \mathcal{P}_{d+1}$ in (12). Then for any $y \in \mathbb{R}^d$ and $z := (y', 1)'$, we have*

$$(16) \qquad u_{\nu',d+1}(z'A^{-1}z) \equiv u_{\nu,d}((y - \mu)'\Sigma^{-1}(y - \mu)).$$

*In particular, this holds for $A = A(\nu')$ and its corresponding $\mu = \mu_\nu \in \mathbb{R}^d$ and $\Sigma = \Sigma_\nu \in \mathcal{P}_d$.*

(d)

$$(17) \qquad \int u_{\nu,d}((y - \mu_\nu)'\Sigma_\nu^{-1}(y - \mu_\nu))dP(y) = 1.$$

(e) *For $h$ defined by (14) with $\rho = \rho_{\nu,d}$, $(\mu_\nu, \Sigma_\nu)$ is the M-functional $\theta_1$ for $P$.*

(f) *If, on the other hand, $P \notin \mathcal{V}_{d,\nu+d}$, then $(\mu, \Sigma) \mapsto Ph(\mu, \Sigma)$ for $h$ as in (e) has no critical points.*

*Proof.* (a): Theorem 8(b) applies since $u_{\nu',d+1}$ satisfies its hypotheses, with $a_0(u_{\nu',d+1}) = \nu' + d + 1 = \nu + d > d + 1$.

(b): By (10), multiplying by $A(\nu')^{-1}$ and taking the trace gives

$$d + 1 = \int u_{\nu',d+1}\left(z'A(\nu')^{-1}z\right)z'A(\nu')^{-1}z dQ(z).$$

We also have, since $z_{d+1} \equiv 1$, that $A(\nu')_{d+1,d+1} = \int u_{\nu',d+1}(z'A(\nu')^{-1}z)dQ(z)$. For any $s \geq 0$, we have $su_{\nu',d+1}(s) + \nu'u_{\nu',d+1}(s) = \nu + d$. Combining gives

$$d + 1 = \nu + d - \nu'\int u_{\nu',d+1}\left(z'A(\nu')^{-1}z\right)dQ(z),$$

and (b) follows.

(c): We can just apply (13) with $\gamma = 1$, and for $A = A(\nu')$, part (b).

(d): This follows from (b) and (c).

(e): By Proposition 10, for $\gamma = 1$ fixed, the relation (12) is a homeomorphism between $\{A \in \mathcal{P}_{d+1} : A_{d+1,d+1} = 1\}$ and $\{(\mu, \Sigma) : \mu \in \mathbb{R}^d, \Sigma \in \mathcal{P}_d\}$. So this also follows from Theorem 8.

(f): We have $\nu + d > d + 1$, so $Q \notin \mathcal{U}_{d+1,\nu+d}$ by Proposition 11. By Theorem 8(a), $Qh$ defined by (8) for $\rho = \rho_{\nu',d+1}$ has no critical point $A$. Suppose $Ph$ has a critical point $(\mu, \Sigma)$ for $\rho = \rho_{\nu,d}$. Let $A := A(\Sigma, \mu, 1) \in \mathcal{P}_{d+1}$. By an affine transformation we can assume $\mu = 0$ and $\Sigma = I_d$, the $d \times d$ identity matrix, so $A = I_{d+1}$. Equations for $\Sigma = I_d$ to be a critical point can be written in the form $\partial/\partial(\Sigma^{-1})_{ij} = 0$, $1 \leq i \leq j \leq d$. By (16) it follows easily that equation (10) holds for $B = A$ and $u = u_{\nu',d+1}$ with the possible exception of the $(d + 1, d + 1)$ entry. Summing the equations for the diagonal $(i, i)$ entries for $i = 1, \ldots, d$, it follows that the $(d + 1, d + 1)$ equation and so (10) holds. By Proposition 6, we get that $A$ is a critical point of the given $Qh$, a contradiction. □

Kent, Tyler and Vardi [16, Theorem 3.1] show that if $u(s) \geq 0$, $u(0) < +\infty$, $u(\cdot)$ is continuous and nonincreasing for $s \geq 0$, and $su(s)$ is nondecreasing for $s \geq 0$, up to $a_0 > d$, and if (17) holds with $u$ in place of $u_{\nu,d}$ at each critical point $(\mu, \Sigma)$ of $Qh$, then $u$ must be of the form $u(s) = u_{\nu,d}(s) = (\nu + d)/(\nu + s)$ for some $\nu > 0$. Thus, the method of relating pure scatter functionals in $\mathbb{R}^{d+1}$ to location-scatter functionals in $\mathbb{R}^d$ given by Theorem 12 for $t$ functionals defined by functions $u_{\nu,d}$ does not extend directly to any other functions $u$.

When $d = 1$, $P \in \mathcal{V}_{1,\nu+1}$ requires that $P(\{x\}) < \nu/(1+\nu)$ for each point $x$. Then $\Sigma$ reduces to a number $\sigma^2$ with $\sigma > 0$. If $\nu > 1$, and $P \notin \mathcal{V}_{1,\nu+1}$, then for some unique $x$, $P(\{x\}) \geq \nu/(\nu+1)$. One can extend $(\mu_\nu, \sigma_\nu)$ by setting $\mu_\nu(P) := x$ and $\sigma_\nu(P) := 0$, with $(\mu_\nu, \sigma_\nu)$ then being weakly continuous at all $P$ [6, Section 6].

The $t_\nu$ functionals $(\mu_\nu, \Sigma_\nu)$ defined in this section can't have a weakly continuous extension to all laws for $d \geq 2$, because such an extension of $\mu_\nu$ would give a weakly continuous affinely equivariant location functional defined for all laws, which is impossible by Theorem 4(b). Here is an example showing that for $d = 2$ and empirical laws with $n = 6$, invariant under interchanging $x$ and $-x$, and/or $y$ and $-y$, so that an affinely equivariant $\mu$ must be 0, there is no continuous extension of the scatter matrix $\Sigma_\nu$ to laws concentrated in lines. For $k = 1, 2, \dots$ let

$$P^{(k)} := \frac{1}{6}\left[\delta_{(-1,-1/k)} + \delta_{(-1,1/k)} + \delta_{(1,-1/k)} + \delta_{(1,1/k)}\right] + \frac{1}{3}\delta_{(0,0)},$$

$$Q^{(k)} := \frac{1}{6}\left[2\delta_{(-1,0)} + \delta_{(0,-1/k)} + \delta_{(0,1/k)} + 2\delta_{(1,0)}\right].$$

Then for each $\nu > 1$, all members of both sequences have mass $\leq 1/3 < \nu/(\nu+2)$ at each point and mass $\leq 2/3 < (\nu+1)/(\nu+2)$ on each line, so are in $\mathcal{U}_{2,\nu+2}$ and the functionals $\mu_\nu, \Sigma_\nu$ are defined for them. By the symmetries $x \leftrightarrow -x$ and $y \leftrightarrow -y$, $\mu_\nu \equiv 0$ and $\Sigma_\nu$ is diagonal on all these laws. Both sequences converge to the limit $P = \frac{1}{3}\left[\delta_{(-1,0)} + \delta_{(0,0)} + \delta_{(1,0)}\right]$, which is concentrated in a line and so is not in $\mathcal{U}_{2,\nu+2}$ for any $\nu$. $\Sigma_\nu(P^{(k)})$ converges to $\binom{a^{(\nu)}\ 0}{0\ \ 0}$ but $\Sigma_\nu(Q^{(k)})$ converges to $\binom{b^{(\nu)}\ 0}{0\ \ 0}$ where $a(\nu) := 2(1 - \nu^{-1})/3 \neq b(\nu) := (2 + \nu^{-1})/3$. We also have $\Sigma_\nu(Q^{(1)}) = \binom{b^{(\nu)}\ 0}{0\ \ c(\nu)}$ with $c(\nu) = \frac{1}{3}(1 - \nu^{-1})$, so that, in contrast to Theorem 5(a), $\Sigma_\nu$ is not proportional to the covariance matrix $\binom{2/3\ \ 0}{0\ \ 1/3}$ for any $\nu < \infty$, but $\Sigma_\nu$ converges to the covariance as $\nu \to +\infty$, as is not surprising since the $t_\nu$ distribution converges to a normal one.

## 5. Differentiability of $t$ functionals

Let $(S, e)$ be any separable metric space, in our case $\mathbb{R}^d$ with its usual Euclidean metric. Recall the space $BL(S, e)$ of all bounded Lipschitz real functions on $S$, with its norm $\|f\|_{BL}$. The dual Banach space $BL^*(S, e)$ has the dual norm $\|\phi\|_{BL}^*$, which metrizes the weak topology on probability measures [5, Theorem 11.3.3].

Let $V$ be an open set in a Euclidean space $\mathbb{R}^d$. For $k = 1, 2, \dots$, let $C_b^k(V)$ be the space of all real-valued functions $f$ on $V$ such that all partial derivatives $D^p f$, for $D^p := \partial^{[p]}/\partial x_1^{p_1} \cdots \partial x_d^{p_d}$ and $0 \leq [p] := p_1 + \cdots + p_d \leq k$, are continuous and bounded on $V$. On $C_b^k(V)$ we have the norm

(18)     $$\|f\|_{k,V} := \sum_{0 \leq [p] \leq k} \|D^p f\|_{\sup,V}, \quad \text{where} \quad \|g\|_{\sup,V} := \sup_{x \in V} |g(x)|.$$

Then $(C_b^k(V), \|\cdot\|_{k,V})$ is a Banach space. For $k = 1$ and $V$ convex in $\mathbb{R}^d$ it's straightforward that $C_b^1(V)$ is a subspace of $BL(V, e)$, with the same norm for $d = 1$ and an equivalent one if $d > 1$.

Substituting $\rho_{\nu,d}$ from (15) into (6) gives for $y \in \mathbb{R}^d$ and $A \in \mathcal{P}_d$,

$$L_{\nu,d}(y, A) := \frac{1}{2} \log \det A + \frac{\nu + d}{2} \log \left[1 + \nu^{-1} y' A^{-1} y\right],$$

so that in (7) we get

$$h_\nu(y, A) := h_{\nu,d}(y, A) := L_{\nu,d}(y, A) - L_{\nu,d}(y, I).$$

Differentiating with respect to entries $C_{ij}$ where $C = A^{-1}$, and recalling $u_{\nu,d}(s) \equiv (\nu + d)/(\nu + s)$, we get as shown in [6, Section 5]

$$(19) \qquad \frac{\partial h_{\nu,d}(y, A)}{\partial C_{ij}} = \frac{\partial L_{\nu,d}(y, A)}{\partial C_{ij}} = -\frac{A_{ij}}{1 + \delta_{ij}} + \frac{(\nu + d) y_i y_j}{(1 + \delta_{ij})(\nu + y' C y)}.$$

For $0 < \delta < 1$ and $d = 1, 2, \ldots$, let

$$\mathcal{W}_\delta := \mathcal{W}_{\delta,d} := \{A \in \mathcal{P}_d : \max(\|A\|, \|A^{-1}\|) < 1/\delta\}.$$

The following is proved in [6, Section 5].

**Lemma 13.** *For any $\delta \in (0, 1)$ let $U := U_\delta := \mathbb{R}^d \times \mathcal{W}_{\delta,d}$. Let $A \in \mathcal{P}_d$ be parameterized by the entries $C_{kl}$ of $C = A^{-1}$. For any $\nu \geq 1$, the functions $\partial L_{\nu,d}/\partial C_{kl}$ in (19) are in $C_b^j(U_\delta)$ for all $j = 1, 2, \ldots$.*

To treat $t$ functionals of location and scatter in any dimension $p$ we will need functionals of pure scatter in dimension $p + 1$, so in the following lemma we only need dimension $d \geq 2$. The next lemma, proved in [6, Section 5], helps to show differentiability of $t$ functionals via implicit function theorems, as it implies that the derivative of the gradient (the Hessian) of $Qh$ is non-singular at a critical point. Let $T(d) := \{(i, j) : 1 \leq i \leq j \leq d\}$.

**Lemma 14.** *For each $\nu > 0$, $d = 2, 3, \ldots$, and $Q \in \mathcal{U}_{d,\nu+d}$, let $A(\nu) = A_\nu(Q) \in \mathcal{P}_d$ be the unique critical point of $Qh(\cdot)$ defined by (8) for $\rho = \rho_{\nu,d}$ defined by (15). For $C = A^{-1}$, the Hessian matrix $\partial^2 Qh(A)/\partial C_{ij} \partial C_{kl}$ with rows indexed by $(i, j) \in T(d)$ and columns by $(k, l) \in T(d)$ is positive definite at $A = A(\nu)$.*

For any $\nu > 0$ and $A \in \mathcal{P}_d$, let $L_{i,j,\nu}(x, A) := \partial L_{\nu,d}(x, A)/\partial C_{ij}$ from (19). Let $X := BL^*(\mathbb{R}^d, e)$ for the usual metric $e(s, t) := |s - t|$. Again, parameterize $A \in \mathcal{P}_d$ with inverse $C$ by $\{C_{ij}\}_{1 \leq i \leq j \leq d}$ in $\mathbb{R}^{d(d+1)/2}$. Consider the open set $\Theta := \mathcal{P}_d \subset \mathbb{R}^{d(d+1)/2}$ and the function $F := F_\nu$ from $X \times \Theta$ into $\mathbb{R}^{d(d+1)/2}$ defined by

$$(20) \qquad F(\phi, A) := \{\phi(L_{i,j,\nu}(\cdot, A))\}_{1 \leq i \leq j \leq d}.$$

Then $F$ is well-defined because $L_{i,j,\nu}(\cdot, A)$ are all bounded and Lipschitz functions of $x$ for each $A \in \Theta$; in fact, they are $C^1$ with bounded derivatives equal except possibly for signs to second partials of $L_{\nu,d}$ with respect to $C_{ij}$. The next fact, proved in [6, Section 5], uses some basic notions and facts from infinite-dimensional calculus, given in [2] and reviewed in the Appendix of [6].

**Theorem 15.** *Let $X := BL^*(\mathbb{R}^d, e)$. In parts (a) through (c), let $\nu > 0$.*

*(a) The function $F = F_\nu$ is $C^\infty$ (for Fréchet differentiability) from $X \times \Theta$ into $\mathbb{R}^{d(d+1)/2}$.*

*(b) Let $Q \in \mathcal{U}_{d,\nu+d}$, and take the corresponding $\phi_Q \in X$. At $A_\nu(Q)$, the $d(d+1)/2 \times d(d+1)/2$ matrix $\partial F(\phi_Q, A)/\partial C := \{\partial F(\phi_Q, A)/\partial C_{kl}\}_{1 \leq k \leq l \leq d}$ is invertible.*

*(c) The functional $Q \mapsto A_\nu(Q)$ is $C^\infty$ for the $BL^*$ norm on $\mathcal{U}_{d,\nu+d}$.*

*(d) For each $\nu > 1$, the functional $P \mapsto (\mu_\nu, \Sigma_\nu)(P)$ given by Theorems 8 and 12 is $C^\infty$ on $\mathcal{V}_{d,\nu+d}$ for the $BL^*$ norm.*

To prove asymptotic normality of $\sqrt{n}(T(P_n) - T(P))$ for $T = (\mu_\nu, \Sigma_\nu)$, the dual-bounded-Lipschitz norm $\|\cdot\|_{BL}^*$ is too strong for some heavy-tailed distributions. Giné and Zinn [12] proved that for $d = 1$, $\{f : \|f\|_{BL} \leq 1\}$ is a $P$-Donsker class if and only if $\sum_{j=1}^\infty \Pr(j - 1 < |X| \leq j)^{1/2} < \infty$ for $X$ with distribution $P$. To define norms better suited to present purposes, for $\delta > 0$ and $r = 1, 2, \ldots$, let $\mathcal{F}_{\delta,r}$ be the set of all functions of $y$ appearing in (19) and their partial derivatives with respect to $C_{ij}$ through order $r$, for any $A \in \mathcal{W}_\delta$. Then each $\mathcal{F}_{\delta,r}$ is a uniformly bounded VC major class as in Theorem 3. Let $Y_{\delta,r}$ be the linear span of $\mathcal{F}_{\delta,r}$. Let $X_{\delta,r}$ be the set of all real-valued linear functionals $\phi$ on $Y_{\delta,r}$ for which $\|\phi\|_{\delta,r} := \sup\{|\phi(f)| : f \in \mathcal{F}_{\delta,r}\} < \infty$. For $A \in \mathcal{W}_{\delta,d}$ and $\phi \in X_{\delta,r}$, define $F(\phi, A)$ again by (20), which makes sense since each $L_{i,j,\nu}(\cdot, A) \in \mathcal{F}_{\delta,r}$ for any $r = 0, 1, 2, \ldots$ by definition.

The next two theorems are also proved in [6, Section 5]. Theorem 17 is a delta-method fact.

**Theorem 16.** *Let $0 < \delta < 1$. For any positive integers $d$ and $r$, Theorem 15 holds for $X = X_{\delta,r+3}$ in place of $BL^*(\mathbb{R}^d, e)$, $\mathcal{W}_{\delta,d}$ in place of $\Theta$, and $C^r$ in place of $C^\infty$ wherever it appears (parts (a), (c), and (d)).*

**Theorem 17.** *(a) For any $d = 2, 3, \ldots$ and $\nu > 0$, let $Q \in \mathcal{U}_{d,\nu+d}$. Then the empirical measures $Q_n \in \mathcal{U}_{d,\nu+d}$ with probability $\to 1$ as $n \to \infty$ and $\sqrt{n}(A_\nu(Q_n) - A_\nu(Q))$ converges in distribution to a normal distribution with mean 0 on $\mathbb{R}^{d(d+1)/2}$ if $A \in \mathcal{P}_d$ is parameterized by $\{A_{ij}\}_{1 \leq i \leq j \leq d}$, or a different normal distribution for the parameterization by $\{A_{ij}^{-1}\}_{1 \leq i \leq j \leq d}$ as above. The limit distributions can also be taken on $\mathbb{R}^{d^2}$, concentrated on symmetric matrices.*

*(b) Let $d = 1, 2, \ldots$ and $1 < \nu < \infty$. For any $P \in \mathcal{V}_{d,\nu+d}$, the empirical measures $P_n \in \mathcal{V}_{d,\nu+d}$ with probability $\to 1$ as $n \to \infty$ and the functionals $\mu_\nu$ and $\Sigma_\nu$ are such that as $n \to \infty$, $\sqrt{n}[(\mu_\nu, \Sigma_\nu)(P_n) - (\mu_\nu, \Sigma_\nu)(P)]$ converges in distribution to some normal distribution with mean 0 on $\mathbb{R}^d \times \mathbb{R}^{d^2}$, whose marginal on $\mathbb{R}^{d^2}$ is concentrated on $d \times d$ symmetric matrices.*

Now, here is a statement on uniformity as $P$ and $Q$ vary, proved in [6, Section 5].

**Proposition 18.** *For any $\delta > 0$ and $M < \infty$, the rate of convergence to normality in Theorem 17(a) is uniform over the set $\mathcal{Q} := \mathcal{Q}(\delta, M)$ of all $Q \in \mathcal{U}_{d,\nu+d}$ such that $A_\nu(Q) \in \mathcal{W}_\delta$ and*

$$Q(\{z : |z| > M\}) \leq (1 - \delta)/(\nu + d), \tag{21}$$

*or in part (b), over all $P \in \mathcal{V}_{d,\nu+d}$ such that $\Sigma_\nu(P) \in \mathcal{W}_\delta$ and (21) holds for $P$ in place of $Q$.*

## Acknowledgements

I thank Lutz Dümbgen and David Tyler very much for kindly sending me copies of their papers [9], [22], and [10]. I also thank Zuoqin Wang and Xiongjiu Liao for helpful literature searches.

\* The author has seen the Rousseeuw (1985) paper cited in secondary sources (MathSciNet and several found by JSTOR) but not in the original.